\documentclass{amsart}

\newif\ifdraft\draftfalse
\newif\ifcite\citefalse
\newif\ifblow\blowtrue

\makeatletter
\@namedef{subjclassname@2020}{\textup{2020} Mathematics Subject Classification}
\makeatother

\usepackage{amsfonts,amssymb, amsmath}
\usepackage{oldgerm}
\usepackage{url}
\usepackage{amscd}
\usepackage{float}
\usepackage{mathrsfs} 
\usepackage[hyperindex]{hyperref}  
\usepackage{tikz}
\usepackage{graphicx}
\usepackage[alphabetic]{amsrefs}
\ifdraft\ifcite\usepackage{showkeys}\else\usepackage[notcite,notref]{showkeys}\fi\fi
\usepackage{mathtools}

\DeclarePairedDelimiter\floor{\lfloor}{\rfloor}

\newtheorem{proposition}[equation]{Proposition}
\newtheorem{theorem}[equation]{Theorem}

\newtheorem{lemma}[equation]{Lemma}

\theoremstyle{definition}
\newtheorem{definition}[equation]{Definition}

\newtheorem{notation}[equation]{Notation}

\theoremstyle{remark}

\newtheorem{remark}[equation]{Remark}

\newtheorem{rem}[equation]{Remark}
\newtheorem{question}[equation]{Question}

\newtheorem{example}[equation]{Example}

\numberwithin{equation}{section}

\def\bc{\begin{cases}}
\def\ec{\end{cases}}

\def\Bbb{\mathbb}

\def\cs{{\mathcal S}}

\def\bc{{\mathbb C}}

\def\bh{{\mathbb H}}

\def\bp{{\mathbb P}}

\def\bz{{\mathbb Z}}

\def\er{\eqref}

\def\bz{\mathbb Z}

\def\bc{\mathbb C}
\def\bp{\mathbb P}
\def\bh{\mathbb H}

\def\lp2{L_pH_{2p}}
\def\bean{\begin{eqnarray}}
\def\eean{\end{eqnarray}}
\def\bea{\begin{eqnarray*}}
\def\eea{\end{eqnarray*}}
\def\beq{\begin{equation}}
\def\eeq{\end{equation}}
\def\beq*{\begin{equation*}}
\def\eeq*{\end{equation*}}
\def\bal{\begin{align*}}
\def\eal{\end{align*}}
\def\baln{\begin{align}}
\def\ealn{\end{align}}

\def\beg{\begin{gather*}}
\def\eng{\end{gather*}}
\def\bqu{\begin{question}}
\def\equ{\end{question}}

\def\ban{\begin{proof}[Answer]}
\def\ean{\end{proof}}

\def\on{\operatorname}

\def\bqu{\begin{question}}
\def\equ{\end{question}}

\def\0110{\begin{matrix} 0 & 1\\1&0\end{matrix}}

\def\ban{\begin{proof}[Answer]}
\def\ean{\end{proof}}
\def\wt{\widetilde}

\def\ben{\begin{equation}}
\def\een{\end{equation}}

\def\j1{{(j+1)}}

\def\f32{{}_3F_2}

\newcommand{\zn}{{\mathbb Z}/N{\mathbb Z}}
\newcommand{\pzn}{{\mathbb P}^1({\mathbb Z}/N{\mathbb Z})}

\begin{document}

\title
{Connected fundamental domains for congruence subgroups}

\author{Zhaohu Nie}
\email{zhaohu.nie@usu.edu}
\address{Department of Mathematics and Statistics, Utah State University, Logan, UT 84322-3900, USA}

\author{C. Xavier Parent}
\email{xavier.parent@usu.edu}
\address{Department of Mathematics and Statistics, Utah State University, Logan, UT 84322-3900, USA}

\subjclass[2020]{11F06, 20H05}

\begin{abstract}
We give explicit sets of right coset representatives for the congruence subgroups $\Gamma_0(N)$, $\Gamma_1(N)$ and $\Gamma(N)$, and prove that the corresponding unions of standard modular triangles are connected fundamental domains.  The construction is based on a study of 
the projective line ${\mathbb P}^1(\zn)$. For every residue class $j\in\zn$, the number of representatives above $j$ is governed by the simple function
\[W_j=\min\{m\in\bz_{>0}\mid mj-1\in (\zn)^*\}.\]
We also include examples illustrating how the connected domains make cusp and boundary data visible.
\end{abstract}

\maketitle

\section{Introduction}
Let $\Gamma(1)=SL_2(\bz)$ act on the upper-half plane $\bh = \{z\,|\,\on{Im}z>0\}$ by the M\"obius transformation
$$
\gamma \cdot z = \frac{az + b}{cz+d},\quad
\gamma=\begin{pmatrix}
a & b\\
c & d
\end{pmatrix}
\in \Gamma(1).
$$

Two points in $\bh$ are called \emph{equivalent} under $\Gamma<\Gamma(1)$ if they are in the same orbit of $\Gamma$. We have the following definition from \cite{Shimura}*{p.~15}.
\begin{definition}\label{Shm} For a discrete subgroup $\Gamma<\Gamma(1)$, we call a set $F\subset \bh$ a \emph{fundamental domain} for $\Gamma$ if 
\begin{enumerate}
\item $F$ is open and connected,
\item no two points of $F$ are equivalent under $\Gamma$,
\item every point of $\bh$ is equivalent to a point in $\overline F$, the closure of $F$.
\end{enumerate}
\end{definition}

It is well known that a fundamental domain for $\Gamma(1)$ in the above sense is 
\begin{equation}\label{first D}
D = \Big\{z\in \bh\,\Big|\, |z|>1, |\on{Re}z|<\frac{1}{2}\Big\}.
\end{equation}
Since $-I$ is in the center of $\Gamma(1)$, we see that $(\pm I) \Gamma$ is also a subgroup of $\Gamma(1)$. When $-I\in \Gamma$, $(\pm I) \Gamma=\Gamma$. 
Note that the M\"obius transformation for $-I$ is the identity. If we have a choice of right coset representatives
\begin{equation}\label{reps}
    (\pm I)\Gamma\backslash \Gamma(1) = \{(\pm I)\Gamma \gamma_1,\dots, (\pm I)\Gamma\gamma_n\},
\end{equation}
where $n=[\Gamma(1):(\pm I)\Gamma]$, we will show in Lemma \ref{fund dom} that 
\begin{equation}\label{intcls}
\on{Int}\Big(\bigcup_{i=1}^n\overline{\gamma_i D} \Big),
\end{equation}
is a fundamental domain for $\Gamma$ provided that it is {\bf connected}, where the closure is, as always in this paper, taken in $\bh$. We will then call such a list of right coset representatives 
$
\{\gamma_1,\dots,\gamma_n\}
$
\emph{fundamental}. 

Consider the congruence subgroups \cite{DS}*{p.~13}
$$
\Gamma(N) < \Gamma_1(N) < \Gamma_0(N) < \Gamma(1),\quad N>1.
$$
The main purpose of this paper is to provide explicit and canonical lists of fundamental right coset representatives such that the corresponding domains are connected.

There are several classical and computational approaches to fundamental domains for congruence subgroups.  Chuman gave an explicit Schreier system of representatives for $\Gamma_0(N)\backslash \Gamma(1)$, with applications to generators and relations of $\Gamma_0(N)$ \cite{Chuman}.  More generally, Kulkarni's Farey-symbol method gives an arithmetic-geometric framework for constructing special polygons for finite-index subgroups of the modular group \cite{Kulkarni}.  There are also computer programs for finding and drawing connected fundamental domains, including Verrill's program \cite{Verrill}.  
The emphasis here is different.  We give a closed-form family of
representatives whose connectedness is built into the construction.
The resulting domains provide a concrete entry point to the geometry
of modular curves and complement the cusp-class and Farey-symbol
approaches of Chuman and Kulkarni.

It is also useful to view our domains at a slightly coarser level.  Each standard triangle has three distinguished sides, which in the figures may be called the left side, the right side, and the base side.  If $\gamma D$ is one of our standard triangles, then the triangle across its base side is $\gamma S D$.  Thus, whenever both $\gamma$ and $\gamma S$ occur in our representative list, the common base side is an internal edge and the two triangles form a two-triangle cell, or kite.  In this sense the domains constructed below may be regarded first as triangulated domains and then, after internal base-side gluings, as partially decomposed into kites.  This gives a visual bridge between our explicit standard-tile domains and the coarser polygonal regions appearing in Farey-symbol constructions.  We do not claim that this provides a new construction of Kulkarni special polygons, as the Farey symbols have become well established and easy to use \cite{KL}. 

In this paper, we mainly work in the ring $\bz/N\bz$. We write the elements as if they are just integers as $a, b\in \bz/N\bz$. Note that the greatest common divisors $\gcd(a, N)$ and $\gcd(a, b, N)$ are well defined as positive integers between 1 and $N$.

To make our final pictures more centered and hence more visualizable, we define our symmetric choice of residue classes mod $N$. The usual choice of classes from $0$ to $N-1$ mathematically would work the same (with suitable small modifications), but the pictures would not be as nice when drawn out. 

\begin{notation}\label{mymod} 
    For $N>1$, let $N_1=\floor{\frac{N-1}{2}}$ and $N_2=\floor{\frac{N}{2}}$, where $\floor{\cdot}$ is the usual floor function. 
    
    For $x\in \bz $, we define $\wt{x}$ to be the unique integer such that $\wt{x}\equiv x \mod N$ and $-N_1\leq \wt{x}\leq N_2$.
\end{notation}

To formulate the construction, we study the projective line over the finite ring
\[
\pzn=\{(a,b)\in (\zn)^2\mid \gcd(a,b,N)=1\}/\sim,
\]
where $(a,b)\sim(a',b')$ if $(a,b)=u(a',b')$ for some unit
$u\in(\zn)^*$.  We write the class of $(a,b)$ as $(a:b)$.

For $x=(a:b)\in\pzn$, define
\begin{equation}\label{eqn:M(a,b)}
M(x)=M(a:b)=\min\{m\in\bz_{\geq 0}\mid ma-b\in(\zn)^*\}.
\end{equation}
The condition $ma-b\in(\zn)^*$ says that $(a:b)$ is distant from the affine point $(1:m)$ in the usual terminology for projective lines over rings; see, for example, \cite{BH}*{\S 2}.  This observation is only a useful way to recognize the condition and will not be used later.

\begin{proposition}\label{prop:M-finite}
The function $M$ in \er{eqn:M(a,b)} is well defined and finite on $\pzn$.
\end{proposition}

Let $m=M(a:b)$ and put
\begin{equation}\label{def c}
c=ma-b\in(\zn)^*.
\end{equation}
Then
\(
(a:b)=(ac^{-1}:bc^{-1})
\)
and
\begin{equation}\label{get 01}
m(ac^{-1})-(bc^{-1})=1\quad\text{in }\zn.
\end{equation}
This gives a distinguished section of the natural projection from admissible pairs to $\pzn$:
\begin{equation}\label{sigma-section}
\sigma:\pzn\longrightarrow(\zn)^2,
\qquad
\sigma(a:b)=c^{-1}(a,b).
\end{equation}
This is clearly independent of the choices of $a, b$ in $(a:b)$, and we call $\sigma(a:b)$ the $M$-normalized representative of $(a:b)$.
If
\(
\sigma(x)=(j,\ell),
\)
then
\begin{equation}\label{def jl}
M(x)j-\ell=1\quad\text{in }\zn.
\end{equation}
Thus every $M$-normalized representative with first coordinate $j$ has the form
$$(j,M(x)j-1).$$

We introduce some notation:
\begin{gather*}
\cs = \sigma(\pzn)\subset (\zn)^2,\\
\cs(j) = \cs\cap \pi_1^{-1}(\{j\}),
\end{gather*}
where $\pi_1:(\zn)^2\to\zn$  is projection to the first component. Therefore $\cs(j)$ is the fiber of $\pi_1:\cs\to \zn$ over $j$, and we aim to characterize it completely now.

\begin{definition}\label{def:W}
For $j\in\zn$, define
\begin{equation}\label{def w}
W_j=\min\{r\in\bz_{>0}\mid rj-1\in(\zn)^*\}.
\end{equation}
\end{definition}
The set in \eqref{def w} is nonempty, since $Nj-1=-1$ in $\zn$.
We regard $W_j$ as the first break in the vertical string
\[
(j,-1),\quad (j,j-1),\quad (j,2j-1),\quad\ldots .
\]
The following proposition says that $\cs(j)$  is exactly the initial segment of this string before the break.

\begin{proposition}\label{prop:initial-segment}
For each $j\in\zn$, we have 
\[\cs(j) =
\{(j,mj-1)\mid 0\leq m<W_j\}.
\]  
\end{proposition}





We define
\[
M_j=W_j-1.
\]
Thus $M_j$ is the last height in the normalized vertical string over $j$, while $W_j$ is the next height, where the string stops.  

We now translate these projective-line representatives into matrices.  Recall that $\Gamma(1)$ is generated by
\begin{equation}\label{ts}
T=\begin{pmatrix}1&1\\0&1\end{pmatrix},
\qquad
S=\begin{pmatrix}0&-1\\1&0\end{pmatrix},
\end{equation}
with $S^2=-I$.

The main result is the following theorem, where a fundamental list for \(\Gamma\) means,
as in \eqref{reps}, a set of representatives for
\((\pm I)\Gamma\backslash \Gamma(1)\).  
\begin{theorem}\label{main thm}
Let $\Gamma(N)<\Gamma_1(N)<\Gamma_0(N)<\Gamma(1)$.
\begin{enumerate}
\item For $\Gamma_0(N)$, a fundamental list of right coset representatives is $\Theta_0$ given by
\begin{equation}\label{oureps}
\begin{cases}
ST^i,& -N_1\leq i\leq N_2,\\
ST^jST^m,& -N_1\leq j\leq N_2,\ \gcd(j,N)>1,\ 0\leq m\leq M_j.
\end{cases}
\end{equation}

\item For $\Gamma_1(N)$, a fundamental list of right coset representatives is $\Theta_1$ given by
\begin{equation}\label{now1}
\begin{cases}
ST^i,& -N_1\leq i\leq N_2,\\
ST^jST^m,& -N_1\leq j\leq N_2,\ \gcd(j,N)>1,\ 0\leq m\leq M_j,\\
ST^kST^i,& -N_1\leq k\leq -2,\ \gcd(k,N)=1,\ -N_1\leq i\leq N_2,\\
ST^kST^{\wt{k^{-1}+j}}ST^m,& -N_1\leq k\leq -2,\ \gcd(k,N)=1,\\
&\qquad -N_1\leq j\leq N_2,\ \gcd(j,N)>1,\ 0\leq m\leq M_j,
\end{cases}
\end{equation}
where $k^{-1}$ denotes an inverse of $k$ modulo $N$, and $\wt{x}$ is defined in Notation \ref{mymod}.

\item For $\Gamma(N)$, a fundamental list of right coset representatives is
\[
\{T^\ell\gamma\mid -N_1\leq\ell\leq N_2,\ \gamma\in\Theta_1\}.
\]
\end{enumerate}
\end{theorem}

Part (1), for $\Gamma_0(N)$, is the main point.  The representatives for $\Gamma_1(N)$ and $\Gamma(N)$ are then obtained by the standard quotient maps between congruence subgroups.

In Section 2, we provide proofs for the results mentioned above. 
In Section 3, we give some examples and pictures.

\section{Proofs of the results}

\begin{proof}[Proof of Proposition \ref{prop:M-finite}]
If $(a,b)$ is replaced by $u(a,b)$ with $u\in(\zn)^*$, then $ma-b$ is replaced by $u(ma-b)$.  Thus the condition $ma-b\in(\zn)^*$ is independent of the chosen representative of $(a:b)$, so $M$ is well defined.

It remains to show that the set in \eqref{eqn:M(a,b)} is nonempty.  
Let $p$ be a prime divisor of $N$.  If $p\mid a$, then $p\nmid b$ by $\gcd(a, b, N)=1$, so $ma-b\not\equiv0\pmod p$ for every $m$.  If $p\nmid a$, then the condition $ma-b\not\equiv0\pmod p$ excludes exactly one residue class of $m$ modulo $p$.  By the Chinese remainder theorem, choosing $m$ to avoid these finitely many bad residue classes modulo the primes dividing $N$, we get $ma-b$ prime to $N$.  Hence $M$ is finite.
\end{proof}

\begin{proof}[Proof of Proposition \ref{prop:initial-segment}]
For $(j,\ell)\in \cs(j)$, we denote $M(j:\ell)=m$. Then by \er{def jl}, $\ell = mj-1$. 

On the other hand, for $m\in \bz_{\geq 0}$, consider $(j, mj-1)\in (\zn)^2$. Suppose that $\forall\, 0\leq t<m,\ tj - (mj-1)\notin (\zn)^*$, then by 
$$
m\cdot j - (mj-1) = 1,
$$
we know from \er{def w} and \er{def jl}, that
$$
M(j: mj-1) = m,\quad\text{and}\quad \sigma(j: mj-1)=(j, mj-1).
$$

Therefore, we have the following equivalence
\begin{align*}
&\negthinspace\qquad\quad(j, mj-1)\in \cs(j)\text{ with }M(j:mj-1) = m\\
&\iff\, \forall 0\leq t<m,\ tj - (mj-1)\notin (\zn)^*\\
&\iff\, \forall 1\leq r\leq m,\ rj-1\notin (\zn)^*\\
&\iff m<W_j.
\end{align*}
\end{proof}

We now record two elementary facts that convert the algebraic representative lists into fundamental domains.

\begin{lemma}\label{fund dom}
The set
\[
\on{Int}\Big(\bigcup_{i=1}^n\overline{\gamma_iD}\Big)
\]
described in \eqref{intcls} is a fundamental domain for \(\Gamma\), provided that it is connected.
\end{lemma}

\begin{proof}
Write
\begin{equation}\label{ED}
E=\on{Int}\Big(\bigcup_{i=1}^n\overline{\gamma_iD}\Big),
\qquad
\mathcal D=\bigcup_{i=1}^n\gamma_iD .
\end{equation}
Since the tiles are ideal geodesic triangles,
$$\overline E=\bigcup_{i=1}^n\overline{\gamma_iD}=\overline{\mathcal D},$$ 
{so every nonempty open subset of }$E$\text{ meets }$\mathcal D$.

We check the conditions in Definition \ref{Shm}. For (1), openness and connectedness of $E$ are immediate from the hypothesis.  For (2), suppose \(x,y\in E\) and \(y=\gamma x\) for some \(\gamma\in\Gamma\). 
Let \(U=E\cap \gamma^{-1}E\), a nonempty open set since it contains $x$. Since
\(\mathcal D\) is dense in \(E\), and \(\gamma^{-1}\mathcal D\)
is dense in \(\gamma^{-1}E\), we may choose
\[
x'\in U\cap \mathcal D\cap \gamma^{-1}\mathcal D .
\]
Then \(y'=\gamma x'\in\mathcal D\). By the definition \er{ED}, write
\(x'=\gamma_i x''\) and \(y'=\gamma_j y''\), with
\(x'',y''\in D\). 
Since \(D\) is a fundamental domain for
\(\Gamma(1)\), we get \(x''=y''\) and
\(\gamma\gamma_i=\pm\gamma_j\).
As the \(\gamma_i\) are representatives for \((\pm I)\Gamma\backslash\Gamma(1)\), this gives \(i=j\), \(\gamma=\pm I\), and hence \(x=y\).

To prove (3), for any \(x\in\bh\), choose \(\theta\in\Gamma(1)\) and \(y\in\overline D\) with \(x=\theta y\).  Writing \(\theta=\pm\gamma\gamma_i\) with \(\gamma\in\Gamma\) by \er{reps}, we get \(x=\gamma(\gamma_i y)\) and \(\gamma_i y\in\overline E\).  Thus every \(\Gamma\)-orbit meets \(\overline E\).
\end{proof}

\begin{definition}\label{G_theta}
For a finite set \(\Theta\subset PSL_2(\bz)=SL_2(\bz)/(\pm I)\), let \(G_\Theta\) be the graph with vertex set \(\Theta\), where \(\theta\) and \(\theta'\) are adjacent if
\[
\theta'=\theta S,
\qquad
\theta'=\theta T,
\qquad\text{or}\qquad
\theta'=\theta T^{-1}.
\]
\end{definition}

\begin{lemma}\label{Graph Connection}
If \(G_\Theta\) is connected, then
\[
\on{Int}\Big(\bigcup_{\theta\in\Theta}\overline{\theta D}\Big)
\]
is connected.
\end{lemma}

\begin{proof}
If two vertices of \(G_\Theta\) are adjacent, the corresponding closed triangles share a full edge: after applying the inverse of one vertex, this is just the fact that \(\overline D\) shares its three edges with \(\overline{SD}\), \(\overline{TD}\), and \(\overline{T^{-1}D}\).  Choose a spanning tree for \(G_\Theta\).  Starting from one tile and adding the remaining tiles along this tree, each new tile is attached along a full edge, so the interior of the union remains connected.  Hence the full tile union has connected interior.
\end{proof}

\begin{proof}[Proof of Theorem \ref{main thm}]
\emph{Part (1), \(\Gamma_0(N)\).}
From \cite{Crem}*{\S 2.2}, the right cosets in \(\Gamma_0(N)\backslash\Gamma(1)\) are in bijection with \(\pzn\).  The bijection is induced by the bottom-row map
\begin{equation}\label{the map}
{\mathcal R}:\Gamma(1)\longrightarrow\zn\times\zn,
\qquad
\begin{pmatrix}a&b\\ c&d\end{pmatrix}\longmapsto(c,d),
\end{equation}
followed by passage to the projective class.

Let
\begin{gather*}
A=\{(a:b)\in\pzn\mid a\in(\zn)^*\}=\{(1:r)\mid r\in\zn\},\\
H=\pzn\setminus A.
\end{gather*}
Define
\begin{equation}\label{pref}
\on{pr}:\pzn\longrightarrow\zn\times\zn,
\qquad
\on{pr}(x)=
\begin{cases}
(1,r), &x=(1:r)\in A,\\
\sigma(x), &x\in H.
\end{cases}
\end{equation}
Thus \(\operatorname{pr}\) keeps the usual affine representative on \(A\) and uses the \(M\)-normalized representative on \(H\) as in \er{sigma-section}. 

We compute
\[
ST^i=\begin{pmatrix}0&-1\\1&i\end{pmatrix}
\]
and
\begin{equation}\label{sts}
ST^jST^m=\begin{pmatrix}-1&-m\\ j&mj-1\end{pmatrix}.
\end{equation}
Therefore
\[
{\mathcal R}(ST^i)=(1,i),
\qquad
{\mathcal R}(ST^jST^m)=(j,mj-1).
\]
The first family represents the affine part \(A\).  The global fibers of \(\sigma\) are described by Proposition \ref{prop:initial-segment} for every residue class \(j\).  
The auxiliary map \(\operatorname{pr}\) deliberately treats the affine part differently: it keeps these points in the usual affine coordinate, for the purpose of keeping the corresponding domain connected.  The remaining non-affine points are exactly the fibers with nonunit first coordinate, and the second family with \(0\leq m\leq M_j\) represents them in their \(M\)-normalized form.  Therefore, the bottom rows of the matrices in \(\Theta_0\) are exactly the pairs \(\operatorname{pr}(x)\) as \(x\) ranges through \(\pzn\).  

It remains to prove connectedness.  In \(G_{\Theta_0}\), the vertex \(S=ST^0\) belongs to \(\Theta_0\), and the vertices \(ST^i\) form the path
\[
(ST^{-N_1},\ldots,ST^{-1},S,ST,\ldots,ST^{N_2}).
\]
For each nonunit \(j\), the vertex \(ST^jS\) is adjacent to \(ST^j\), and the path
\[
(ST^jS,ST^jST,\ldots,ST^jST^{M_j})
\]
connects every vertex above \(j\) to \(ST^jS\).  Hence every vertex of \(G_{\Theta_0}\) is connected to \(S\).  Lemmas \ref{Graph Connection} and \ref{fund dom} show that \(\Theta_0\) is a fundamental list for \(\Gamma_0(N)\).

\medskip
\emph{Part (2), \(\Gamma_1(N)\).}
We use the exact sequence \cite{DS}*{p.~14}
\begin{align*}
1\to (\pm I)\Gamma_1(N)\to \Gamma_0(N)&\to (\zn)^*/(\pm I)\to 1,\\
\begin{pmatrix}a&b\\ c&d\end{pmatrix}&\mapsto \{\pm d\}\pmod N.
\end{align*}
For a unit \(k\in\zn\), choose an inverse \(k^{-1}\) modulo \(N\).  Then
\[
ST^kST^{k^{-1}}S=
\begin{pmatrix}-k^{-1}&1\\ k^{-1}k-1&-k\end{pmatrix}
\in\Gamma_0(N),
\]
and, as \(k\) runs through representatives for \((\zn)^*/(\pm I)\), these matrices give representatives for \((\pm I)\Gamma_1(N)\backslash\Gamma_0(N)\).  For the class of \(k=1\), we use \(I\) instead of \(STSTS=-T^{-1}\).  Thus
\begin{equation}\label{rep10}
\{I,ST^kST^{k^{-1}}S\mid -N_1\leq k\leq -2,\ \gcd(k,N)=1\}
\end{equation}
is a set of representatives for \((\pm I)\Gamma_1(N)\backslash\Gamma_0(N)\).

Multiplying these representatives by the elements of \(\Theta_0\) gives the list \(\Theta_1\) in \eqref{now1}.  The first two lines are \(I\cdot\Theta_0\). The third and the fourth lines follow from
\begin{align*}
(ST^kST^{k^{-1}}S)(ST^i)&=-ST^kST^{k^{-1}+i},\\
(ST^kST^{k^{-1}}S)(ST^jST^m)
   &=-ST^kST^{k^{-1}+j}ST^m,
\end{align*}
where we replace the exponents $k^{-1}+i$ and \(k^{-1}+j\) by their symmetric representative,
as in Notation \ref{mymod}. We can do this, since the right cosets  represented by $ST^i$ and $ST^jST^m$ in $\Theta_0$
\er{oureps} depend only on the residue classes of $i$ and $j$ modulo $N$.

The graph \(G_{\Theta_1}\) contains \(G_{\Theta_0}\).  Since the allowed \(k\)'s lie in the range of the \(i\)'s, \(ST^kS\) is adjacent to \(ST^k\), and
\[
(ST^kST^{-N_1},\ldots,ST^kS,\ldots,ST^kST^{N_2})
\]
connects every vertex \(ST^kST^i\) to \(S\).  Similarly, \(\wt{k^{-1}+j}\) lies in the same range, and
\[
(ST^kST^{\wt{k^{-1}+j}}S,ST^kST^{\wt{k^{-1}+j}}ST,\ldots,
ST^kST^{\wt{k^{-1}+j}}ST^{M_j})
\]
connects every remaining vertex to \(S\).  Hence \(G_{\Theta_1}\) is connected, and Lemmas \ref{Graph Connection} and \ref{fund dom} show that \(\Theta_1\) is a fundamental list for \(\Gamma_1(N)\).

\medskip
\emph{Part (3), \(\Gamma(N)\).}
The matrices \(T^\ell\), with \(-N_1\leq\ell\leq N_2\), give representatives for \(\Gamma(N)\backslash\Gamma_1(N)\).  Therefore the products \(T^\ell\gamma\), \(\gamma\in\Theta_1\), give representatives for \((\pm I)\Gamma(N)\backslash\Gamma(1)\).  Since \(I=ST^0ST^0\) belongs to \(\Theta_1\) in \(PSL_2(\bz)\), the vertices \(T^\ell\) occur in the list and form the path
\[
(T^{-N_1},\ldots,T^{-1},I,T,\ldots,T^{N_2}).
\]
For each fixed \(\ell\), left multiplication by \(T^\ell\) gives a copy of the connected graph \(G_{\Theta_1}\) attached at \(T^\ell\).  Hence the full graph is connected.  Lemmas \ref{Graph Connection} and \ref{fund dom} finish the proof.
\end{proof}

\begin{remark}\label{pq}
     We consider a special case when $N=p^rq^s$ has two prime factors $p$ and $q$. For $\gcd(a, N)>1$, $\gcd(b, N)>1$ and $\gcd(a, b, N)=1$, we show that $M(a:b)=1$, that is, $a-b$ is a unit mod $N$. This is true, since under our assumption, $\gcd(a, N)$ is one of the prime factors to some positive power and $\gcd(b, N)$ is the other. Therefore $a-b$ is not divisible by either of the prime factors, so a unit mod $N$. Hence, $M$ can be at most 1 in this case. 

    On the other hand, if $N=p^r$ is a prime power, then $\gcd(a, N)>1$ and $\gcd(a, b, N)=1$ would require $\gcd(b, N)=1$, and then $M(a:b)=0$. 
    
    Meanwhile, the values $M_j$ are not bounded by $1$ in general; the first example with larger values is illustrated for $N=30$ in Example \ref{ex30} below.
\end{remark}

\section{Examples}

We implemented the domains in Maple by drawing the ideal geodesic triangles.  In the figures, Maple colors edges red, blue and green, with overlap priority red$\,>\,$blue$\,>\,$green.

\begin{example}\label{ex6}
We show a labeled picture of our fundamental domain for $\Gamma_0(6)$.  Since
\[
[\Gamma_0(N):(\pm I)\Gamma_1(N)]=\frac{\phi(N)}2
\]
and $\phi(6)/2=1$, this is also the domain for $\Gamma_1(6)$.

\begin{figure}[ht]
\begin{center}
\includegraphics[width=0.95\textwidth,height=0.68\textheight,keepaspectratio]{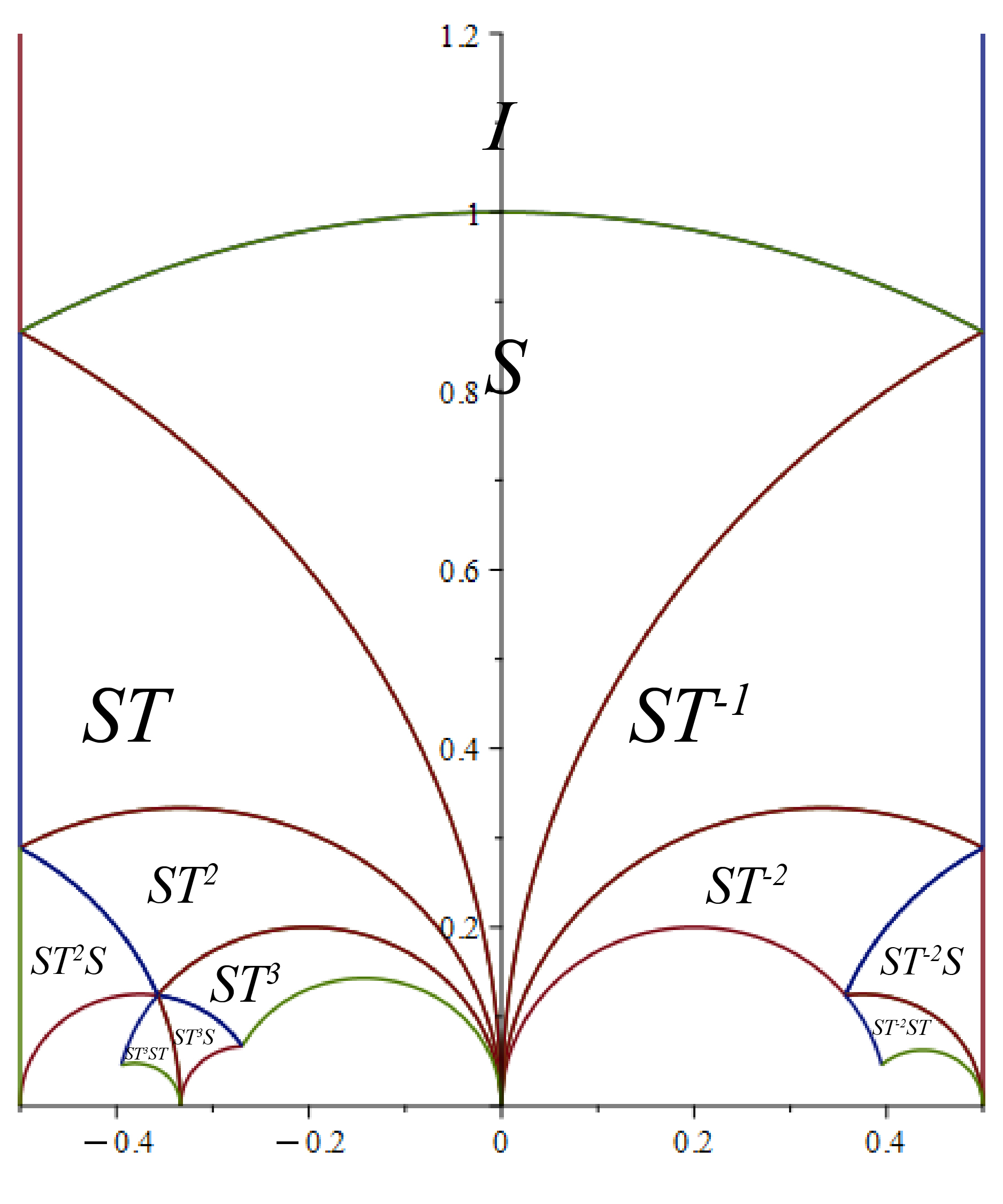}
\caption{Our Fundamental Domain for $\Gamma_0(6)$ Labeled}\label{fig:gamma06}
\end{center}
\end{figure}

The functions $M_j$ and $W_j$ are defined for every residue class $j\in\bz/6\bz$.  The nonunit columns, which are the columns entering the second family of representatives in Theorem \ref{main thm}, are
\begin{center}
\begin{tabular}{|c|c|c|c|c|}
\hline
$j$ & $-2$ & 0 & 2 & 3 \\
\hline
$M_j$ & 1 & 0 & 0 & 1\\
\hline
$W_j$ & 2 & 1 & 1 & 2\\
\hline
\end{tabular}
\end{center}
The unit columns are part of the same globally defined functions, but the corresponding points are represented through the affine part of the list.

The connected picture makes several features visible at once.  Internal base-side pairs $(ST^jD,ST^jSD)$ form kites whenever both representatives occur; for instance, the $j=0$ pair is $(SD,D)$ because $S^2=-I$ acts trivially on $\bh$.  Apart from the pairing of the left and right vertical sides, the bottom boundary arcs labeled $1,\ldots,8$ are paired by
\[
1\sim 8,
\qquad 2\sim 7,
\qquad 3\sim 4,
\qquad 5\sim 6,
\]
as in \cite{N}*{Theorem 1.21}.  The cusps are represented by $\infty,0,1/2,1/3$, with widths $1,6,3,2$; there are no elliptic points of orders $2$ or $3$, and the genus formula gives $g(X_0(6))=0$.
\end{example}

\begin{example}\label{ex30}
The case $N=30=2\cdot3\cdot5$ is the first one in which values $M_j>1$ occur.  Here
\[
[\Gamma(1):\Gamma_0(30)]=\psi(30)=30\prod_{p\mid 30}\left(1+\frac1p\right)=72.
\]
The affine part of $\bp^1(\bz/30\bz)$ has $30$ elements.  The classes $(j:-1)$ with nonunit $j$ contribute $22$ more, and the remaining $20$ non-affine classes account for the higher values of $M_j$.

For nonunit $j\in\bz/30\bz$, the values are as follows.  These are again only the nonunit columns of the globally defined functions $M_j$ and $W_j$.
\begin{center}
\begin{tabular}{|c|l|c|}
\hline
$M_j$ & $j$ & $W_j$\\
\hline
0 & $-12, -10, -6, 0, 2, 8, 12, 14$ & 1\\
\hline
1 & $-14, -9, -8, -5, -3, 4, 6, 9, 10, 15$ & 2\\
\hline
2 & $-4, -2$& 3\\
\hline
3 & 3, 5&4\\
\hline
\end{tabular}
\end{center}
Thus the number of representatives of the form $ST^jST^m$ with $j$ nonunit and $0\leq m\leq M_j$ is
\[
1\cdot 8+2\cdot 10+3\cdot 2+4\cdot 2=42=\psi(30)-30.
\]

For $ST^jST^m$ in \er{sts}, the associated cusp is
\[
(ST^jST^m)(\infty)=-\frac1j,
\]
and it occurs $M_j+1$ times.  Grouping the above values by cusp equivalence gives the following widths; the contributing $j$'s are counted with multiplicity $M_j+1$.
\begin{center}
\begin{tabular}{|c|c|l|}
\hline
cusp rep & width & contributing $j$'s\\
\hline
$\frac12$ & 15 & $2,4,8,14,-14,-8,-4,-2$\\
$\frac13$ & 10 & $3,9,-9,-3$\\
$\frac15$ & 6 & $5,-5$\\
$\frac16$ & 5 & $6,12,-12,-6$\\
$\frac1{10}$ & 3 & $10,-10$\\
$\frac1{15}$ & 2 & $15$\\
$\infty$ & 1 & affine cusp at $\infty$\\
$0$ & 30 & affine cusp at $0$\\
\hline
\end{tabular}
\end{center}
These widths agree with the standard computations in \cite{DS}*{Prop. 3.8.3} and \cite{S}*{Algorithm 1.19}.  This example is meant only to record the first occurrence of the higher values $M_j>1$; the general pattern is proved in the companion paper \cite{N}.
\end{example}

\medskip

\end{document}